\documentclass{article}

\usepackage{arxiv}

\usepackage[utf8]{inputenc} 
\usepackage[T1]{fontenc}    
\usepackage{hyperref}       
\usepackage{url}            
\usepackage{booktabs}       
\usepackage{amsfonts}       
\usepackage{nicefrac}       
\usepackage{microtype}      
\usepackage{lipsum}

\usepackage[final]{graphicx}

\title{THE HILBERT SPACE BASIS AND HILBERT'S EIGHTH PROBLEM
}

\author{
 Kirill V.~Kapitonets \\
 BAUMAN MSTU, GRADUATE 1990\\
  MCC EuroChem\\
  Moscow\\
  Russian Federation\\
  \texttt{kkapitonets@live.com} \\
}

\begin{document}
\maketitle
\begin{abstract}
The paper considers the Hilbert space $\hat{H}_r$ of real functions summable with the square $L^2(a,b)_r$ on any interval $\{(a,b)_r\}_{r=1}^{\infty}\in \mathbb{R}$.
\par
It is shown on the basis of the theorem on zeros of real orthogonal polynomials if in $\hat{H}_r$ there exists a complete orthonormal basis $\{f(x)_k\}_{k=1}^{\infty}$ and the function $f(x)\in\{f(x)_k\}_{k=1}^{\infty}$ has zeros, then these zeros are simple and real.
\par
The generalized Hardy function $Z(\sigma,t)=\Re\zeta(\sigma+it)e^{i\theta(t)}$ is considered.
\par
It is shown that in the Hilbert space $\hat{H}_r$ there exists a complete basis $\{Z(\lambda_k,t\}_{k=1}^{\infty}$ where $\lambda_k\in\mathbb{Q}$ and $Z(t)\in\{Z(\lambda_k,t\}_{k=1}^{\infty}$ when $\lambda_k=1/2$, hence the Hardy function $Z(t)=\zeta(1/2+it)e^{i\theta(t)}$ has all simple and real zeros.

\end{abstract}

\section{Introduction}\nopagebreak
In 1859, when Riemann laid the foundation of analytical number theory with his report at the Berlin Academy of Sciences, mathematics in its modern form was just beginning to take shape.
\par
In particular, the axioms of linear space \cite{PE} were formulated by Peano only in 1888.
\par
Therefore, when Riemann discovered a functional dependency between the values of $\zeta(s)$ on two different lines of the complex plane, Riemann does not talk about the linear dependence of these values, although he explicitly points to it.
\par
Using the residue theory of a function of a complex variable, Riemann shows that the residue values (\ref{int}) that are used to calculate the contour integral that Riemann derived to calculate the values of $\zeta(s)$ are proportional, i.e. linearly dependent on the value of the Dirichlet series term that determines the value of $\zeta(s)$ on a parallel line of complex plane.
\begin{equation}\label{int}2\sin\pi s\Pi(s-1)\zeta(s)=2\pi^s\sum n^{1-s}((-i)^{s-1}+i^{s-1})\end{equation}
\par
\begin{figure}[ht!]
\centering
\includegraphics[scale=0.5]{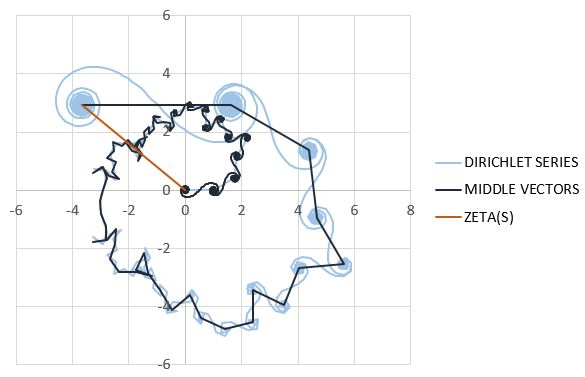}
\caption{The Riemann spiral and the reverse Riemann spiral}
\label{fig:reverse_spiral_approx}
\end{figure}

This linear dependence is also visible in the construction of the Riemann spiral (the segments connect the points corresponding to successive partial sums of the Dirichlet series Fig. \ref{fig:reverse_spiral_approx}), while the middle segments that connect the visible centers of the vortices of the Riemann spiral make up the inverse Riemann spiral and correspond to the values of the residues that Riemann deduced.
\par
Riemann simply had no chance to prove that all zeros of the function of the real variable $\Xi(t)$\footnote{The values of  $\Xi(t)$ determine values of $|\zeta(s)|$ on the critical line, where $\zeta(s)$ on the complex plane does not exist another except for the line itself, on which the values of $\zeta(s)$ will be linearly dependent.}, all are simple and real, which determines that all nontrivial zeros of $\zeta(s)$ lie on this critical line.
\par
Such a chance appeared only in 1932, when Siegel literally dug up in Riemann's posthumous notes \cite{SI} a formula that is now called the Riemann-Siegel formula and is used as the main tool of the theory of $\zeta(s)$.
\par
In particular, using this formula, Karatsuba proved theorems about the density of zeros $\zeta(s)$ on almost all very short intervals of the critical line \cite{KA}.
\par
It is believed that the results obtained by Karatsuba cannot be improved in terms of growth, therefore, unfortunately, the possibilities of analytical number theory are exhausted, because estimates of the density of zeros $\zeta(s)$ become comparable with estimates of the argument $\zeta(s)$ and it is no longer possible to obtain arbitrarily significant the result.
\par
It is also impossible to get a significantly better result of estimating the total number of zeros $\zeta(s)$ on a critical line of more than about 40\%\cite{HU}.
\par
What happened in mathematics by 1932, besides the appearance of the concept of linear dependence?
\par
Firstly, in 1873 Cantor laid the foundations of set theory \cite{CA} and gave mathematicians such a powerful tool that at the beginning of the XX century Bourbaki (a group of French mathematicians) attempted to rewrite the foundations of mathematics \cite{BU}, which has not been completed to this day, because the process of development of mathematics is endless.
\par
Obviously, the concept of infinity as a countable infinity, i.e. the number of natural and rational numbers is Cantor's main achievement.
\par
Cantor showed that the cardinalities of the set of rational and real numbers do not coincide.
\par
Secondly, Dedekind, in 1872, gives a definition of the real number \cite{DE}, which later becomes the final tool for determining the fundamental sequence.
\par
But perhaps the most important event of mathematics in the XX century is the work of Hilbert\cite{HI}, in which in 1912 he laid the foundations of linear functional spaces, which now bear his name - Hilbert spaces.
\par
Hilbert substantiates theoretically trigonometric Fourier series and lays the foundations of generalized Fourier series, which are based on the linear independence of the functions that make up the basis of Hilbert space.
\par
The main result here is the statement about the separability of the Hilbert space, which allows us to assert that with the help of generalized Fourier series in Hilbert space, any function of this space can be represented and such a series will be the Fourier series of this analytic function.
\par
Unfortunately, all the results of Hilbert spaces were very quickly generalized to the field of complex numbers $\mathbb{C}$.
\par
Thus, the chances of proving Riemann's statement again became zero.
\par
In parallel with theoretical research in mathematics, practical calculations are carried out and in 1854 Chebyshev wrote a paper about mechanisms that in their work deviate least from the direct line \cite{CH}, in which he develops the theory of orthogonal polynomials to approximate the expression of functions describing the kinematics of such mechanisms.
\par
Now Chebyshev polynomials are one of the types of classical orthogonal polynomials that approximate functions of a real variable on a given interval $(a,b)\in\mathbb{R}$.
\par
It is important to note that the theorem on zeros of orthogonal polynomials \cite{GA}, which states that these zeros are simple and real, is true only for real polynomials.
\par
Thus, the logic of things tells us that in order to prove the Riemann hypothesis \cite{RI} that all zeros $\Xi(t)$ are simple and real, it is necessary to leave from analytical methods of functions of a complex variable to the field of real numbers $\mathbb{R}$.
\par
\section{The basis of Hilbert space}\nopagebreak
Consider the separable Hilbert space $\hat{H}_r$ of real functions summable with the square $L^2(a,b)_r$ on any interval $\{(a,b)_r\}_{r=1}^{\infty}\in\mathbb{R}$.
\par
Thus, the norm \cite{SA} is defined in $\hat{H}_r$
\begin{equation}\label{norm}||f(x)||=\Big(\int_a^b f(x)^2dx\Big)^{1/2}<\infty\end{equation}
\par
and the scalar product
\begin{equation}\label{prod}(f(x),g(x))=\int_a^b f(x)g(x)dx<\infty\end{equation}
\par
If $\hat{H}_r$ is separable, then there are many complete isomorphic bases in it, then there is a complete universal basis
\begin{equation}\label{base}\{1, x, x^2, x^3{...}x^{m-1}\}_{m=1}^{\infty}\end{equation}
\par
Then there is a complete universal orthogonal basis $\{\phi(x)_n\}_{n=1}^{\infty}$, which can be represented by a general Fourier series
\begin{equation}\label{fourier}\phi(x)_n=\sum_{m=1}^{\infty}a_{nm}x^{m-1}\end{equation}
\par
where
\begin{equation}\label{coeff}a_{nm}=\int_a^b\phi(x)_n x^{m-1}dx\end{equation}
Fourier coefficients of $\phi(x)_n$.
\par
Therefore, for any $\phi(x)_n$, the following inequality will hold
\begin{equation}\label{epsilon}\Big|\phi(x)_n-\sum_{m=1}^{M}a_{nm}x^{m-1}\Big|<\epsilon\end{equation}
\par
Obviously, $\sum_{m=1}^{M}a_{nm}x^{m-1}$ is a polynomial in contrast to the Fourier series (\ref{fourier}).
\par
We use this property of the Hilbert space to prove the statement about the zeros of the basis of the Hilbert space $\hat{H}_r$, based on the theorem about zeros of \cite{GA} real orthogonal polynomials.
\par
STATEMENT 1\nopagebreak
\par
Let $\{f(x)_k\}_{k=1}^{\infty}$ be an arbitrary complete orthogonal basis $\hat{H}_r$ on any interval $\{(a,b)_r\}_{r=1}^{\infty}\in\mathbb{R}$ different from the complete universal orthogonal basis $\{\phi(x)_n\}_{n=1}^{\infty}$.
\par
Then if function $f(x)\in\{f(x)_k\}_{k=1}^{\infty}$ has zeros, then these zeros are simple and real.
\par
PROOF\nopagebreak
\par
If $\{f(x)_k\}_{k=1}^{\infty}$ is a complete orthogonal basis of $\hat{H}_r$ on any interval $\{(a,b)_r\}_{r=1}^{\infty}\in\mathbb{R}$ other than the full universal orthogonal basis $\{\phi(x)_n\}_{n=1}^{\infty}$, then each function $f(x)_k$ can be represented by a Fourier series in the full universal orthogonal basis $\{\phi(x)_n\}_{n=1}^{\infty}$ on any interval $\{(a,b)_r\}_{r=1}^{\infty}\in\mathbb{R}$
\begin{equation}\label{fourier_phi}f(x)_k=\sum_{n=1}^{\infty}b_{nk}\phi(x)_n\end{equation}
\par
where
\begin{equation}\label{coeff_phi}b_{nk}=\int_a^b f(x)_n\phi(x)_ndx\end{equation}
Fourier coefficients of $f(x)_k$.
\par
Then, obviously, the following inequality
\begin{equation}\label{epsilon_phi}\Big|f(x)_k-\sum_{n=1}^{N}b_{nk}\phi(x)_n\Big|<\epsilon\end{equation}
\par
Combine the inequalities (\ref{epsilon}) and (\ref{epsilon_phi})
\begin{equation}\label{epsilon_phi_2}\Big|f(x)_k-\sum_{n=1}^{N}b_{nk}\sum_{m=1}^{M}a_{nm}x^{m-1}\Big|<\epsilon\end{equation}
\par
Obviously, $\sum_{n=1}^{N}b_{nk}\sum_{m=1}^{M}a_{nm}x^{m-1}$ is a polynomial in contrast to the Fourier series (\ref{fourier_phi}).
\par
Let's arbitrarily choose two functions $f_1,f_2\in\{f(x)_k\}_{k=1}^{\infty}$.
\par
If on any interval $\{(a,b)_r\}_{r=1}^{\infty}\in\mathbb{R}$
\begin{equation}\label{polinomals}p_1=\sum_{n=1}^{N}b_{n1}\sum_{m=1}^{M}a_{nm}x^{m-1},  p_2=\sum_{n=1}^{N}b_{n2}\sum_{m=1}^{M}a_{nm}x^{m-1}\end{equation}
\par
Therefore, \nopagebreak
\begin{equation}\label{epsilons}|f_1-p_1|<\epsilon, |f_2-p_2|<\epsilon\end{equation}
\par
Then, obviously, at any interval $\{(a,b)_r\}_{r=1}^{\infty}\in\mathbb{R}$ $p1\perp p2$ at $N\rightarrow\infty, M\rightarrow\infty$, because
\begin{equation}\label{polinomals_prod}(p_1,p_2)\rightarrow\int_a^b f_1f_2dx=0, N\rightarrow\infty,  M\rightarrow\infty\end{equation}
\par
Let $f_1$ have zeros $\{\alpha_k\}$ on the interval $(a,b)\in\{(a,b)_r\}_{r=1}^{\infty}$ , then $p_1$ also has zeros $\{\beta_k\}$ and
\begin{equation}\label{zeros}\beta_k\rightarrow\alpha_k, N\rightarrow\infty,  M\rightarrow\infty\end{equation}
\par
Then by the theorem on zeros of real orthogonal polynomials $\{\beta_k\}$ are simple and real, hence $\{\alpha_k\}$ are also simple and real.
\par
Statement 1 is proved.
\par
\section{Generalized Hardy function}\nopagebreak
First, let us formulate the main statement that we will prove.
\par
STATEMENT 2\nopagebreak
\par
In the Hilbert space $\hat{H}_r$ on any interval $\{(a,b)_r\}_{r=1}^{\infty}\in\mathbb{R}$ there is a complete basis $\{Z(\lambda_k,t\}_{k=1}^{\infty}$ where $\lambda_k\in\mathbb{Q}$, such that the function $Z(1/2,t)\in\{Z(\lambda_k,t\}_{k=1}^{\infty}$, which defines $|\zeta(s)|$ on the critical line, has all simple and real zeros.
\par
It is obvious that $Z(1/2,t)=\zeta(1/2+it)e^{i\theta(t)}$ is a Hardy function that Siegel returned from oblivion in 1932, and Karatsuba called it \cite{KAR,GAB}.
\begin{equation}\label{theta}\theta(t)=\frac{t}{2}\log(\frac{t}{2\pi})-\frac{t}{2}-\frac{\pi}{8}+\frac{1}{48t}+\frac{7}{5760t^3}+\frac{31}{80640t^5}+O(t^{-7})\end{equation}
\par
We obtain the formula for the Hardy function in trigonometric form directly from the Dirichlet series.
\begin{equation}\label{zeta}\zeta(s)=\sum_{n=1}^{\infty}n^{-s}=\sum_{n=1}^{\infty}n^{-\sigma}(\cos(t\log n)-i\sin(t\log n))\end{equation}
\begin{equation}\label{exp_theta}e^{i\theta(t)}=\cos \theta(t)+i\sin \theta(t)\end{equation}
Multiply the two functions and apply the trigonometric formulas for the difference of the angles of cosines and sines.
\begin{equation}\label{zeta_theta}\zeta(s)e^{i\theta(t)}=\sum_{n=1}^{\infty}n^{-\sigma}\cos(\theta(t)-t\log n)+i\sum_{n=1}^{\infty}n^{-\sigma}\sin(\theta(t)-t\log n)\end{equation}
\par
At $\sigma=1/2$\nopagebreak
\begin{equation}\label{sin_theta}\sum_{n=1}^{\infty}n^{-1/2}\sin(\theta(t)-t\log n)\equiv 0\end{equation}
Then
\begin{equation}\label{cos_theta}\zeta(1/2+it)e^{i\theta(t)}=\sum_{n=1}^{\infty}n^{-1/2}\cos(\theta(t)-t\log n)\end{equation}
\par
Formally, the series (\ref{cos_theta}) diverges, because the series (\ref{zeta}) diverges at the value $\sigma=1/2$, but $\zeta(s)$ is defined on the entire complex plane, so the methods of summing divergent series are applicable to the series (\ref{cos_theta}).
\par
Including from the series (\ref{cos_theta}) Riemann obtained the now well-known Riemann-Siegel formula
\begin{equation}\label{z_theta}Z(t)=2\sum_{n=1}^{N}n^{-1/2}\cos(\theta(t)-t\log n)+R(t), N=\Big\lfloor\sqrt{\frac{t}{2\pi}}\Big\rfloor\end{equation}
\par
\begin{figure}[ht!]
\centering
\includegraphics[scale=0.5]{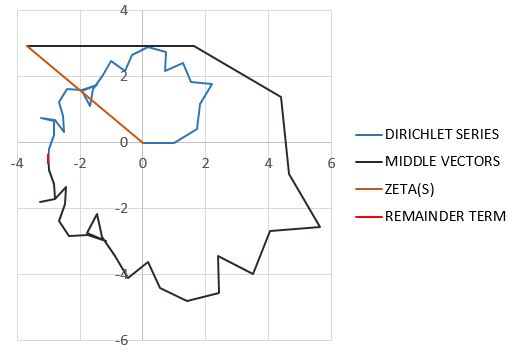}
\caption{A finite system of vectors of an approximate functional equation}
\label{fig:finite_vector_system}
\end{figure}

Like the functional equation of the Riemann Zeta function, the Riemann-Siegel formula has a geometric representation Fig\ref{fig:finite_vector_system}.
\par
The non-intersecting parts of the Riemann spiral and the inverse Riemann spiral are connected by a segment corresponding to the residual term $R(t)$.
\par
The intermediate transformation (\ref{zeta_theta}) can be viewed in two relative ways:\nopagebreak
\par
1) rotation of the vector $\zeta(s)$ by the angle $\theta(t)$;
\par
2) rotation of the complex plane by an angle $-\theta(t)$.
\par
Then in the second case, the real part of the expression (\ref{zeta_theta}) will denote the projection of the vector $\zeta(s)$ on the direction of the unit vector $e^{i\theta(t)}$, and the imaginary part - on the perpendicular direction, i.e. on the direction of the unit vector $ie^{i\theta(t)}$.
\par
Obviously, each projection of the vector $\zeta(s)$ can be considered as a separate real function.
\par
Hence, we can generalize the Hardy function to the entire complex plane.
\par
DEFINITION 1\nopagebreak
\par
Let's denote the real part of the product $\zeta(s)e^{i\theta(t)}$ for a fixed value $\sigma$ by the sign $Z(\sigma, t)$ and call this real function the generalized Hardy function.
\par
Then
\begin{equation}\label{cos_theta_sigma}Z(\sigma, t)=\Re\zeta(\sigma+it)e^{i\theta(t)}=\sum_{n=1}^{\infty}n^{-\sigma}\cos(\theta(t)-t\log n)\end{equation}
\par
Formally series\footnote{We defined a generalized Hardy function at a fixed value of $\sigma$, so we got a set of functions $\{Z(\sigma,t)\}$, $\sigma\in\mathbb{R}$.} (\ref{cos_theta_sigma}) diverges at $\sigma\le 1$, because the series (\ref{zeta}) diverges at the value of $\sigma\le 1$, but $\zeta(s)$ is defined on the entire complex plane, therefore, methods for summing divergent series at any values of $\sigma$ are applicable to the series (\ref{cos_theta_sigma}).
\par
Obviously, the Hardy function belongs to the set of generalized Hardy functions for $\sigma=1/2$.
\par
\begin{equation}\label{z_theta_sigma}Z(t)=Z(1/2, t)\end{equation}
STATEMENT 3\nopagebreak
\par
In the Hilbert space $\hat{H}_r$ on any interval $\{(a,b)_r\}_{r=1}^{\infty}\in\mathbb{R}$ there is a complete basis $\{Z(\lambda_k,t\}_{k=1}^{\infty}$, $\lambda_k\in\mathbb{Q}$, such that the function $Z(1/2,t)\in\{Z(\lambda_k,t\}_{k=1}^{\infty}$.
\par
PROOF\nopagebreak
\par
First, we show that $\{Z(\sigma,t)\}_{\sigma\in\mathbb{R}}\in\hat{H}_r$, i.e. on any interval $\{(a,b)_r\}_{r=1}^{\infty}\in\mathbb{R}$ there is a norm
\begin{equation}\label{norm_z}||Z(\sigma,t)||=\Big(\int_a^b Z(\sigma,t)^2dt\Big)^{1/2}<\infty\end{equation}
\par
It follows from (\ref{zeta_theta}) that $\Re\zeta(\sigma+it)e^{i\theta(t)}$ is a harmonic function that, according to the harmonic functions theorem \cite{KART}, up to a constant uniquely defines its holomorphic function $\zeta(\sigma+it)e^{i\theta(t)}$.
\par
$\zeta(s)$ is a entire holomorphic function, therefore, by the mean theorem \cite{KART} is bounded on the boundary of any circle.
\par
Thus, since $\zeta(s)$ is bounded, and $|e^{i\theta(t)}|=1$, then $\Re\zeta(\sigma+it)e^{i\theta(t)}$ is bounded at a fixed $\sigma$ on any interval $\{(a,b)_r\}_{r=1}^{\infty}\in\mathbb{R}$.
\par
In addition, $\zeta(s)$ is defined on the entire complex plane, so the methods of summing divergent series for any values of $\sigma$ are applicable to the series (\ref{cos_theta_sigma}).
\par
Therefore, on any interval $\{(a,b)_r\}_{r=1}^{\infty}\in\mathbb{R}$ for any values of $\sigma$ there is an integral
\begin{equation}\label{z_theta_sigma_int}\int_a^b Z(\sigma,t)^2dt<\infty\end{equation}
\par
Then $\{Z(\sigma, t)\}_{\sigma\in\mathbb{R}}\in\hat{H}_r$.
\par
Now we show that the functions $\{Z(\sigma,t)\}_{\sigma\in\mathbb{R}}$ are linearly independent on any interval $\{(a,b)_r\}_{r=1}^{\infty}\in\mathbb{R}$ for any values of $\sigma_1\ne 1-\sigma_2$.
\par
Since $\zeta(s)$ is defined on the entire complex plane, therefore, the summation methods of divergent series are applicable to the series (\ref{cos_theta_sigma}) for any values of $\sigma$ and, and $|e^{i\theta(t)}|=1$, then we can substitute the series (\ref{zeta_theta}) into the functional equation of the Riemann Zeta function \cite{TI}:
\begin{equation}\label{zeta_theta_func}\sum_{n=1}^{\infty}\frac{\cos(\theta(t)-t\log n)}{n^{\sigma}}+i\sum_{n=1}^{\infty}\frac{\sin(\theta(t)-t\log n)}{n^{\sigma}}=\chi(s)\sum_{n=1}^{\infty}\frac{\cos(\theta(t)-t\log n)}{n^{1-\sigma}}+i\sum_{n=1}^{\infty}\frac{\sin(\theta(t)-t\log n)}{n^{1-\sigma}}\end{equation}
\par
It can be noticed that the functional equation also holds for $\{Z(\sigma,t)\}_{\sigma\in\mathbb{R}}$:
\begin{equation}\label{z_theta_func}\sum_{n=1}^{\infty}n^{-\sigma}\cos(\theta(t)-t\log n)=\chi(s)\sum_{n=1}^{\infty}n^{\sigma-1}\cos(\theta(t)-t\log n)\end{equation}
\par
where\nopagebreak
\begin{equation}\label{chi}\chi(s)=\frac{2^{s-1}\pi^s}{\cos\frac{\pi s}{2}\Gamma(s)}\end{equation}
\par
From (\ref{z_theta_func}), it is obvious that $\{Z(\sigma, t)\}_{\sigma\in\mathbb{R}}$ are linearly dependent at $\sigma_1=1-\sigma_2$ and linearly independent at $\sigma_1\ne 1-\sigma_2$.
\par
Therefore, $\{Z(\lambda_k, t)\}$ are linearly independent for $\lambda_m\ne 1-\lambda_n$, $\lambda_k\in\mathbb{Q}$.
\par
Obviously, we can set by anyway $\{\lambda_k\}_{k=1}^{\infty}$, $\lambda_k\in\mathbb{Q}$ so that $1/2\in\{\lambda_k\}_{k=1}^{\infty}$.
\par
Then $\{Z(\lambda_k,t\}_{k=1}^{\infty}$, $\lambda_k\in\mathbb{Q}$ is the complete basis in $\hat{H}_r$ on any interval $\{(a,b)_r\}_{r=1}^{\infty}\in\mathbb{R}$ and $Z(1/2,t)\in\{Z(\lambda_k,t\}_{k=1}^{\infty}$.
\par
Statement 3 is proved.
\par
\section{Hilbert's eighth problem}\nopagebreak
Hilbert's eighth problem contains two problems that Hilbert attributed to number theory:\nopagebreak
\par
1) the Riemann hypothesis;
\par
2) Goldbach's problem.
\par
Obviously, our paper concerns the Riemann hypothesis, it was just important to draw a parallel between two mathematical concepts associated with the name of one person.
\par
Now that we have practically proved statement 2, it becomes clear that this problem has not been solved for more than 160 years, because:\nopagebreak
\par
1) at the time of the formulation of the hypothesis \cite{RI}, mathematics did not yet have a tool to show that all zeros of the real function $\Xi(t)$ are simple and real;
\par
2) the hypothesis was reformulated by Hilbert in 1900\cite{HIL}, and the formulation was fixed in 2000 as a problem of a function of a complex variable, although it was originally formulated as a problem of a function of a real variable;
\par
3) Lehmer in 1956\cite{LE} discovers the first problematic pair of zeros of the Hardy function and formulates an assumption about the possibility of complex zeros in this real function;
\par
4) In 2003, Ivech\footnote{Aleksandar Ivić  is one of the few authors of a monograph on the Riemann Zeta function.}\cite{IV} finally fixes the phenomenon of Lehmer as one of the main reasons for the refutation of the Riemann hypothesis (although no pair of complex zeros of the Hardy function has been found so far, in other words, no evidence, just Riemann's word against Lehmer's word).
\par
Let's analyze statement 2 again:\nopagebreak
\par
In the Hilbert space $\hat{H}_r$ on any interval $\{(a,b)_r\}_{r=1}^{\infty}\in\mathbb{R}$ there is a complete basis $\{Z(\lambda_k,t\}_{k=1}^{\infty}$, $\lambda_k\in\mathbb{Q}$, such that the function $Z(1/2,t)\in\{Z(\lambda_k,t\}_{k=1}^{\infty}$, which defines $|\zeta(s)|$ on the critical line, has all simple and real zeros.
\par
It was possible to limit ourselves to a simple indication of a real function that has all simple and real zeros.
\par
But such a return to the original formulation cannot in any way affect the possibility of solving this problem, although of the entire content of statement 2, we have only this part left to prove - $Z(1/2,t)$ has all simple and real zeros.
\par
Let 's follow the whole chain of reasoning:\nopagebreak
\par
1) it is impossible to say with complete certainty about a real function that the function has all simple and real zeros, so there is a deep concern about the existence of Lehmer\cite{LE} pairs, where a function of a real variable may have complex roots (or may not have);
\par
2) for functions of a complex variable, there is a more powerful analytical tool \cite{KART}, but nevertheless in this form it is also difficult to solve the problem of zeros $\zeta(s)$;
\par
3) there is a theorem about zeros of real orthogonal polynomials \cite{GA}, but it is also impossible to apply it to analytic functions of a real variable, because they are not polynomials;
\par
4) it is obvious that there is a close relationship between real analytic functions and real orthogonal polynomials, this relationship is so close that orthogonal polynomials best approximate these functions \cite{CH};
\par
5) the most important thing is that this connection between real analytic functions and real orthogonal polynomials manifests itself in normalized linear spaces with a scalar product, but not in simple, but infinite-dimensional \cite{HI,SA}, where the properties of linear functional spaces manifest themselves most clearly, namely, any countable basis of a Hilbert space is complete, that is, it is enough to show that the basis is countable and there is no need to prove that it is complete;
\par
6) but the most important thing in Hilbert spaces is that linear functionals have unconditional continuity, i.e. polynomials continuously \cite{SA} give their properties to analytic functions of Hilbert space, i.e. it is enough to show that the functions belong to Hilbert space and form its basis.
\par
We used clause 6 in proving statement 1, thus we obtained a criterion for substantiating simple and real zeros of a real function, but to apply this criterion, such a function alone is not enough, we need an orthogonal basis, which includes the function under study.
\par
We used the functional equation \cite{TI} and clause 5 to prove the statement 3 about the existence of such a basis, and we were left with one delicate question, to show that the orthogonalization procedure does not affect the property that we extracted with such difficulty, namely, the function of a real variable has simple and real zeros.
\par
Let's prove the last statement.\nopagebreak
\par
STATEMENT 4\nopagebreak
\par
The Gram-Schmidt orthogonalization procedure has no effect on the function belonging to the Hilbert space bases to have simple and real zeros.
\par
PROOF\nopagebreak
\par
The Gram-Schmidt procedure \cite{HO} consists in converting a non-orthogonal basis $\{f(x)_k\}_{k=1}^{\infty}$ into an orthogonal basis $\{g(x)_k\}_{k=1}^{\infty}$, to do this, follow these steps:\nopagebreak
\par
1) take the countable number of linearly independent functions $\{f(x)_k\}_{k=1}^{\infty}$;
\par
2) choose an arbitrary function $f(x)_1\in\{f(x)_k\}_{k=1}^{\infty}$;
\par
3) $g(x)_1=f(x)_1$;
\par
4) $g(x)_i=f(x)_i-\sum_{k=1}^{i-1}PROJ(f(x)_i, g(x)_k)|_{i=2}^{\infty}$, $PROJ(f(x)_i, g(x)_k)$ is the projection of the vector $f(x)_i$ onto the vector $g(x)_k$;
\par
From step \#3 of the Gram-Schmidt procedure, it is obvious that $f(x)_1$ in the process of orthogonalization of the basis $\{f(x)_k\}_{k=1}^{\infty}$ remains unchanged, therefore, we can always choose the function under study first and the orthogonalization procedure it will have no effect on the function belonging to the Hilbert space basis to have simple and real zeros.
\par
Statement 4, and with it statement 2 are proved.
\par
Hence, we can state that the function $Z(1/2,t)\in\{Z(\lambda_k,t\}_{k=1}^{\infty}$, $\lambda_k\in\mathbb{Q}$, which defines $|\zeta(s)|$ on the critical line, has all simple and real zeros, because $\{Z(\lambda_k,t\}_{k=1}^{\infty}$, $\lambda_k\in\mathbb{Q}$ is a basis in the Hilbert space $\hat{H}_r$ and the Gram-Schmidt procedure does not affect the function, belonging to the basis of Hilbert space $\hat{H}_r$, have simple and real zeros.
\par
This means that the Hardy function $Z(t)=Z(1/2, t)$ has simple and real zeros and no Lehmer pairs can become it complex zeros.
\par
\section{Conclusions}\nopagebreak
The Hardy function $Z(t)=Z(1/2, t)$ has simple and real zeros.
\par
Now this is understandable, because it is clear WHY they should be all simple and real.
\par
To answer this question briefly, because at any interval $\{(a,b)_r\}_{r=1}^{\infty}\in\mathbb{R}$ real polynomials in the Hilbert space $\hat{H}_r$ continuously approximate real functions in this Hilbert space
\begin{equation}\label{final_1}\Big|f(x)-\sum_{n=1}^{N}b_{nk}\sum_{m=1}^{M}a_{nm}x^{m-1}\Big|<\epsilon\end{equation}
\par
at $N\rightarrow\infty, M\rightarrow\infty$
\par
thus, first the properties of orthogonal real functions are transferred to real polynomials, which also continuously become orthogonal, that is, if we take any two functions $f_1,f_2\in\{f(x)_k\}_{k=1}^{\infty}$, where $\{f(x)_k\}_{k=1}^{\infty}$ is simple (that is, not orthogonal) the basis of the Hilbert space $\hat{H}_r$, then using the Gram-Schmidt procedure we can find the function $g_2$ orthogonal to the function $f_1$
\begin{equation}\label{final_2}g_2=f_2-PROJ(f_2, f_1)\end{equation}
\par
where $PROJ(f_2, f_1)$ - projection of the vector $f_2$ onto the vector $f_1$
\par
then
\begin{equation}\label{final_3}(p_1,p_2)\rightarrow\int_a^b f_1g_2dx=0, N\rightarrow\infty,  M\rightarrow\infty\end{equation}
\par
where
\begin{equation}\label{final_4}p_1=\sum_{n=1}^{N}b_{n1}\sum_{m=1}^{M}a_{nm}x^{m-1},  p_2=\sum_{n=1}^{N}b_{n2}\sum_{m=1}^{M}a_{nm}x^{m-1}\end{equation}
\par
then
\begin{equation}\label{final_5}|f_1-p_1|<\epsilon, |g_2-p_2|<\epsilon\end{equation}
\par
then orthogonal real polynomials transfer their property to have simple and real zeros to orthogonal real functions, which make up the basis of the Hilbert space $\hat{H}_r$
\begin{equation}\label{final_6}\beta_k\rightarrow\alpha_k, N\rightarrow\infty,  M\rightarrow\infty\end{equation}
\par
where $\{\beta_k\}$ is the orthogonal zeros of the polynomial $p_1$, and $\{\alpha_k\}$ is the zeros of the function $f_1$, then we can conclude that if such functions of a real variable on any interval $\{(a,b)_r\}_{g=1}^{\infty}\in\mathbb{R}$ have zeros, then these zeros are simple and real.
\par
The statement, which was unthinkable in 1859, obviously became possible only thanks to a qualitative change in mathematics, which was achieved thanks to the axiomatic approach (it is strange why this was not done earlier).
\par
Now linear functional spaces allow you to work wonders, which previously mathematicians were selected solely thanks to intuition.
\par
Hilbert insisted on the need to axiomatize mathematics.
\par
Hilbert laid the foundations of linear functional spaces.
\par
Hilbert collected mathematical problems.
\par
What is not clear is why no one has ever asked the question WHY, why zeros should be real or why zeros can be complex.
\par
So now we can say that if $Z(1/2,t)$ belongs to the basis $\{Z(\lambda_k,t\}_{k=1}^{\infty}$, $\lambda_k\in\mathbb{Q}$ of the Hilbert space $\hat{H}_r$, then we can choose any function $Z(\lambda_k,t), \lambda_k\ne 1/2$ belonging to this basis and find a function orthogonal to $Z(1/2,t)$ on any interval $\{(a,b)_r\}_{r=1}^{\infty}\in\mathbb{R}$
\begin{equation}\label{final_7}G(\lambda_k,t)=Z(\lambda_k,t)-PROJ(Z(\lambda_k,t), Z(1/2,t))\end{equation}
\par
Then we get two orthogonal polynomials
\begin{equation}\label{final_8}p_1=\sum_{n=1}^{N}b_{n1}\sum_{m=1}^{M}a_{nm}x^{m-1},  p_2=\sum_{n=1}^{N}b_{n2}\sum_{m=1}^{M}a_{nm}x^{m-1}\end{equation}
\par
Therefore, we can state that $Z(1/2,t)$ and any function $Z(\lambda_k,t))$\footnote{Instead of $Z(1/2,t)$, we can substitute any function $Z(\lambda_k,t)$ in the first place, so we are not talking about the orthogonal function $G(\lambda_k,t)$, but about any function $Z(\lambda_k,t)$.}, $\lambda_k\ne 1/2$, which belongs to the basis $\{Z(\lambda_k,t\}_{k=1}^{\infty}$, $\lambda_k\in\mathbb{Q}$ of the Hilbert space $\hat{H}_r$ if has zeros on any interval $\{(a,b)_r\}_{r=1}^{\infty}\in\mathbb{R}$, then these zeros are simple and real.
\par
Before completing the presentation, let us consider several consequences of the result obtained.
\par
CONSEQUENCE 1\nopagebreak
\par
In 2000, when the Clay Institute awarded a prize for proving the Riemann hypothesis, an extended Riemann hypothesis was added to the problem statement, which concerns Dirichlet L-series.
\par
Obviously, if any Dirichlet character $\{\chi_w\}_{w=1}^{\infty}$ does not violate the linear independence of the modified extended Hardy function $\hat{Z}(\sigma,t,\chi_w)$ if $\sigma_1\ne 1-\sigma_2$, then the extended Riemann hypothesis is also true, because in this case we can also construct the basis $\{\hat{Z}(\lambda_k, t, \chi_w)\}_{k=1}^{\infty}, \lambda_k\in\mathbb{Q}$Hilbert space $\hat{H}_r$, such that $\hat{Z}(1/2,t,\chi_w)\in\{\hat{Z}(\lambda_k, t, \chi_w)\}_{k=1}^{\infty}$, hence all zeros are $\hat{Z}(1/2,t,\chi_w)$ will be simple and real.
\par
CONSEQUENCE 2\nopagebreak
\par
It is known that the function proposed in 1936 by Davenport and Heilbronn \cite{DA} has a Riemann-type functional equation, nevertheless, the Riemann hypothesis does not hold for it, i.e. this function has zeros outside the critical line.
\par
The fact is that the Davenport-Heilbronn function is a linear combination of Dirichlet L-functions, and as Karatsuba\cite{KARA} noted, any linear combination of Dirichlet L-functions has a common multiplier, and, consequently, the values of these functions will be linearly dependent for any values of $\sigma$, then we cannot for such functions, specify a countable basis in the Hilbert space $\hat{H}_r$, which means that the Davenport-Heilbronn function does not have to have all zeros on the critical line.
\par
CONSEQUENCE 3\nopagebreak
\par
In 1975, Voronin proved the theorem on the universality of the Riemann Zeta function \cite{VOR}, which states that functions of a complex variable that have no zeros and are continuous up to the boundary of a circle are approximated by an incomplete Euler product over a finite set of primes.
\par
Obviously, this universality follows from the existence in the Hilbert space of $\hat{H}_r$ on any interval $\{(a,b)_r\}_{r=1}^{\infty}\in\mathbb{R}$ of the complete basis $\{Z(\lambda_k,t\}_{k=1}^{\infty}$ where $\lambda_k\in\mathbb{Q}$, which is defined by the Riemann Zeta function.
\par

\bibliographystyle{unsrt}  


\begin{thebibliography}{1}
\bibitem{HU}Hung Bui, More than 41 percent of the zeros of the zeta function are on the critical line, 1002.4127 arxiv.org
\end{thebibliography}

\end{document}